\documentclass{IEEEtran}
\usepackage{cite}
\usepackage{amsmath,amssymb,amsfonts}
\usepackage{algorithmic}
\usepackage{graphicx}
\usepackage{textcomp}
\usepackage{mathtools}
\usepackage{algorithm}                  
\usepackage{algorithmic}
\usepackage{colortbl}
\usepackage{enumerate}
\usepackage{xcolor}
\def\BibTeX{{\rm B\kern-.05em{\sc i\kern-.025em b}\kern-.08em
    T\kern-.1667em\lower.7ex\hbox{E}\kern-.125emX}}

\newtheorem{theorem}{Theorem}
\newtheorem{lemma}{Lemma}
\newtheorem{proposition}{Proposition}
\newtheorem{assumption}{Assumption}
\newtheorem{definition}{Definition}
\newtheorem{remark}{Remark}
\newtheorem{corollary}{Corollary}
\begin{document}

\title{Stochastic Optimal LQR Control of Discrete-time Systems with Delay and Quadratic Constraints}

\author{Dawei Liu, Juanjuan Xu, and Huanshui Zhang, \IEEEmembership{Fellow, IEEE}
\thanks{This work was supported by the Original Exploratory Program Project of the National Natural Science Foundation of China under Grant 62450004, and by the National Natural Science Foundation of China under Grants 62573262 and 62250056. (Corresponding author: Huanshui Zhang.)}
\thanks{Dawei Liu and Juanjuan Xu are with the School of Control Science and Engineering, Shandong University, Jinan, Shandong 250061, China (e-mail: dwliu@mail.sdu.edu.cn, juanjuanxu@sdu.edu.cn).}
\thanks{Huanshui Zhang is with the College of Electrical Engineering and Automation, Shandong University of Science and Technology, Qingdao 266590, Shandong, China (e-mail: hszhang@sdu.edu.cn).}
\thanks{This work has been submitted to the IEEE for possible publication. Copyright may be transferred without notice, after which this version may no longer be accessible.}}

\IEEEpubid{0000--0000/00\$00.00~\copyright~2021 IEEE}

\maketitle

\begin{abstract}
This article explores discrete-time stochastic linear quadratic regulation (LQR) control with input delay and quadratic constraints.
By using Lagrangian duality, we derive a parameterized optimal controller and cost function via the Riccati-ZXL equation, and determine the optimal parameters through a projected gradient ascent algorithm with global convergence of the dual objective value.
The primary contribution is to formulate the optimal control as a feedback law driven by the state's conditional expectation, where the feedback gain is obtained from the Riccati-ZXL equation associated with the optimal parameters.
Numerical examples demonstrate the effectiveness of the obtained results.
\end{abstract}

\begin{IEEEkeywords}
Stochastic LQR control, quadratic constraints, time delay, projected gradient ascent.
\end{IEEEkeywords}

\section{Introduction}
\IEEEPARstart{L}{inear} quadratic regulation (LQR) is a core method in optimal control theory, widely applied in fields such as engineering, economics, and biology.
Early studies by Bellman \cite{Be:58}, Kalman \cite{Ka:60}, and Letov \cite{Le:61} laid the foundation for LQR control.
As research progressed, extending deterministic LQR control to stochastic systems became crucial in engineering development and practice (see Wonham \cite{Wo:68}, Athans \cite{At:71}, Davis \cite{Da:77}, Bensoussan \cite{Be:92}, and others).
Recent studies have gone further by considering systems influenced by multiplicative noise and time delays, offering explicit solutions for optimal controllers and addressing numerous practical challenges (see \cite{Zh:15}).
For a deeper exploration of these extensions, refer to the works of \cite{Zh:16} and \cite{Zh:17}.
All of the above are unconstrained LQR problems.

Constrained problems arise naturally in optimal control, particularly when multiple performance objectives must be satisfied simultaneously.
Although the classical LQR framework is widely used for its elegant formulation of control design as the minimization of a quadratic cost, real-world applications often involve balancing several objectives at once, and merely minimizing a single cost function typically fails to meet other critical performance requirements.
To explore the efficient frontier of such multi-objective problems, a common approach is to optimize one performance index while imposing constraints on the others.
As a result, a solid understanding of constrained LQR and the development of efficient algorithms become essential, especially for LQR problems with quadratic constraints, since many practical specifications can be naturally expressed in quadratic form. Typical examples include flight planning \cite{St:83}, where fuel-efficient choices of route, altitude, and speed must still ensure timely and safe arrival under strict performance limits, as well as the control of space structures and industrial processes \cite{To:89}, where similar quadratic constraints naturally arise.

Given the significance of LQR problems with quadratic constraints, many researchers have conducted extensive studies in this area. For continuous-time systems, \cite{Li:99} extensively explored optimal LQR control with integral quadratic constraints, and \cite{Su:18} studied optimal LQR control with fixed terminal states and integral quadratic constraints. For stochastic systems, \cite{Zhou:99} focused on stochastic LQR control and studied systems with integral quadratic constraints under indefinite control weights. For discrete-time systems, \cite{Ts:20} studied finite-horizon risk-constrained LQR, while \cite{Yo:23} addressed the infinite-horizon case via a policy optimization approach. Furthermore, \cite{We:25} addressed the nonlinear finite-horizon optimal control problem with isoperimetric constraints by locally linearizing the dynamics and solving a sequence of LQ time-varying subproblems.

All the aforementioned quadratically constrained optimal LQR problems are formulated for delay-free system dynamics, since the incorporation of delays significantly increases analytical complexity.
This is because existing parameter optimization methods are difficult to apply directly to the problem studied in this paper.
For example, \cite{Zhou:99} derived gradient formulas for continuous-time systems via first-order optimality conditions, thereby bypassing the need for an explicit iterative convergence analysis.
\cite{Yo:23} constructed dual subgradients through controller substitution and established convergence based on primal-dual error bounds for a single-constraint control problem.
In \cite{Li:25}, the system under consideration contains no state-dependent multiplicative noise, and the input-delay history prior to the initial time is assumed to be zero.
In this case, the corresponding optimal parameter problem can be reformulated as a semidefinite programming problem (SDP).
In contrast, this paper studies a more general stochastic optimal control problem with input delay, state-dependent multiplicative noise, and multiple quadratic constraints.
This setting gives rise to a coupled discrete-time Riccati-ZXL equation, and consequently neither the SDP method in \cite{Li:25}
\newpage
\noindent
nor the aforementioned gradient or subgradient methods are directly applicable.

To overcome these difficulties, we develop a duality-based framework for the quadratically constrained stochastic LQR problem with input delay. In the finite-horizon case, strong duality is established under Slater's condition. In particular, for the infinite-horizon case, although the primal problem is infinite-dimensional, we establish strong duality under the proposed assumptions. Hence, the primal constrained problem can be equivalently reformulated as a parameterized unconstrained stochastic LQR problem with input delay, together with an associated parameter optimization problem. For each fixed parameter, the parameterized optimal controller and cost function are explicitly constructed from the solution to the Riccati-ZXL equation. The parameter is optimized by a projected gradient ascent algorithm. To prove its convergence, we directly differentiate the parameterized optimal value and establish the boundedness and Lipschitz continuity of both the Riccati-ZXL solution and its parameter derivatives on relevant bounded parameter sets. These estimates yield the Lipschitz continuity of the dual gradient on the bounded level set containing the iterates, and, together with dual concavity, imply the global convergence of the generated dual objective values. The resulting optimal controller is a feedback law driven by the conditional expectation of the state, with the feedback gain determined by the Riccati-ZXL equation associated with the optimal parameter.

The layout of this paper is as follows:
Section II discusses discrete-time stochastic optimal LQR control with delay and quadratic constraints for the finite-horizon case.
Section III addresses solutions for the infinite-horizon scenario, with numerical examples presented in Section IV.
Finally, Section V offers a brief concluding discussion.

\textbf{Notation.}
$\mathbb{R}^n$ denotes the $n$-dimensional Euclidean space;
$\mathbb{R}_+^r$ denotes the $r$-dimensional nonnegative real vector space;
$I$ represents the identity matrix;
The superscript $\top$ indicates the transpose of a matrix;
A symmetric matrix $M \succ 0$ (resp. $M \succeq 0$) means that it is positive definite (resp. positive semidefinite);
$\frac{\partial \varphi(\lambda)}{\partial \lambda_j}$ denotes the partial derivative of $\varphi$ with respect to $\lambda_j$;
$Z_{k,j}\triangleq\frac{\partial Z_k}{\partial \lambda_j}$ denotes the partial derivative of $Z_k$ with respect to $\lambda_j$;
$\big\{\Omega, \mathcal{F}, \mathcal{P}, \{\mathcal{F}_k\}_{k \ge 0}\big\}$ denotes a complete probability space on which a scalar white noise $\omega_k$ is defined such that $\left\{\mathcal{F}_k\right\}_{k \geq 0}$ is the natural filtration generated by $\omega_k$, i.e. $\mathcal{F}_k=\sigma\left\{\omega_0, \ldots, \omega_k\right\}$, augmented by all the $\mathcal{P}$-null sets in $\mathcal{F}$ \cite{Ch:10};
$\hat{x}_{k \mid m} \doteq E\left[x_k \mid \mathcal{F}_{m-1}\right]$ denotes the conditional expectation of $x_k$ with respect to $\mathcal{F}_{m-1}$;
$\|\cdot\|$ denotes the Euclidean norm for vectors and the corresponding induced operator norm for matrices;
$[x]_{+} \doteq \max (x, 0)$ denotes the positive part operator.

\section{Finite-horizon constrained stochastic optimal LQR control}
\subsection{Problem Statement}
Consider the following discrete-time system:
\begin{equation}\label{equ:1-1}
x_{k+1}=\left(A+\omega_k \bar{A}\right) x_k+\left(B+\omega_k \bar{B}\right) u_{k-d},
\end{equation}
where $x_k \in \mathbb{R}^n$ is the state, $u_k \in \mathbb{R}^m$ is the input control with delay $d>0$. $\omega_k$ is a scalar random white noise with zero mean and variance $\sigma^2$. Matrices $A, B, \bar{A}$ and $\bar{B}$ are constant and have compatible dimensions. The initial values $x_0$ and $u_i$, $i=-d,\ldots,-1$, are given.

In the finite-horizon case, the cost function $(i=0)$ and the constraint
functions $(i=1,\ldots,r)$ are defined as
\begin{equation}\label{equ:1-2}
\begin{aligned}
J_i(\mathbf{u})
= \mathbb{E} \Bigg[
& \sum_{k=0}^{N} x_k^\top Q_i x_k
  + \sum_{k=d}^{N} u_{k-d}^\top R_i u_{k-d}  \\
& + x_{N+1}^\top F_i x_{N+1}
\Bigg],
\end{aligned}
\end{equation}
where $\mathbb{E}$ denotes the mathematical expectation with respect to
the noise sequence $\{\omega_0,\omega_1,\ldots\}$, $R_0 \succ 0$,
$R_i \succeq 0$ $(i=1,\ldots,r)$, $Q_i \succeq 0$, and
$F_i \succeq 0$ $(i=0,\ldots,r)$. Moreover, $N$ denotes the horizon length.

Given $c_1, \dots, c_r \in \mathbb{R}$, the constrained stochastic optimal LQR control problem for the finite-horizon case as follows:

\noindent\textbf{Problem 1} Given the dynamical system \eqref{equ:1-1}, find an $\mathcal{F}_{k-1}$-adapted $u_k$ such that $J_0(\mathbf{u})$ is minimized subject to the quadratic constraints $J_i(\mathbf{u}) \leq c_i(i=1, \ldots, r)$.

\subsection{Finite-Horizon Stochastic Optimal Control}
\textbf{Problem 1} is a convex optimization problem. Under the Slater condition, strong duality holds \cite{Lu:97}, enabling the optimal controls to be determined. The Slater condition is defined as follows:
\begin{assumption}[Slater Condition]\label{asp:1-1}
There exists an $\mathcal{F}_{k-1}$-adapted $\bar{u}_k$ such that for all $i = 1, \ldots, r$,
\begin{equation*}
{J}_i(\bar{\mathbf{u}}) < c_i.
\end{equation*}
\end{assumption}

For $\lambda = (\lambda_1,\ldots,\lambda_r)^\top \in \mathbb{R}_+^r$, define the Lagrangian:
\begin{equation*}
J(\mathbf{u},\lambda)
= J_0(\mathbf{u})
  + \sum_{i=1}^{r}\lambda_i\big(J_i(\mathbf{u})-c_i\big),
\end{equation*}
where $\lambda$ is referred to as the Lagrange multiplier. It follows that
\begin{equation*}
\begin{aligned}
J(\mathbf{u},\lambda)
= \mathbb{E}\Big[
& \sum_{k=0}^{N} x_k^\top Q(\lambda)x_k
  + \sum_{k=d}^{N} u_{k-d}^\top R(\lambda)u_{k-d} \\
& + x_{N+1}^\top F(\lambda)x_{N+1}
\Big] - \lambda^\top c ,
\end{aligned}
\end{equation*}
where $c=\left(c_1,\ldots,c_r\right)^\top$, and $Q(\lambda) = Q_0 + \sum\limits_{i=1}^{r} \lambda_i Q_i$, with $R(\lambda)$ and $F(\lambda)$ defined analogously.

Given any $\lambda \in \mathbb{R}_+^r$, we have the associated dual function
\begin{equation}\label{equ:1-3}
\varphi(\lambda)=\min _\mathbf{u} J(\mathbf{u}, \lambda)
\end{equation}

This is a parameterized, unconstrained finite-horizon stochastic optimal LQR control problem with delay and multiplicative noise. The following theorem is a direct result from \cite{Zh:15}.

\begin{theorem}\label{thm:1-1}
For any given $\lambda\in\mathbb{R}_+^r$, $\varphi(\lambda)$ admits a unique optimal control. Moreover, $\Upsilon_k\succ0$ for $k=N,N-1,\ldots,d$. The optimal $\mathcal F_{k-1}$-adapted controller is given by
\begin{equation}\label{equ:1-4}
u_k=-\Upsilon_{k+d}^{-1}M_{k+d}\hat{x}_{k+d\mid k},
\end{equation}
for $k=0,1,\ldots,N-d$, where
\begin{equation*}
\begin{aligned}
\hat{x}_{k+d\mid k}
&\triangleq \mathbb{E}\left[x_{k+d}\mid\mathcal{F}_{k-1}\right]  \\
&= A^d x_k+\sum_{i=1}^{d}A^{i-1}Bu_{k-i}.
\end{aligned}
\end{equation*}

The corresponding dual function value is
\begin{equation}\label{equ:1-5}
\begin{aligned}
\varphi(\lambda)
= \mathbb{E}\Big[
& \sum_{k=0}^{d-1} x_k^\top Q(\lambda)x_k
  + x_d^\top X_d x_d  \\
& - \sum_{i=0}^{d-1}
    x_d^\top (A^\top)^i L_{d+i}A^i
    \hat{x}_{d\mid i}
\Big] - \lambda^\top c .
\end{aligned}
\end{equation}
Here
\begin{equation*}
\begin{aligned}
\hat{x}_{d\mid i}
&\triangleq \mathbb{E}\left[x_d\mid\mathcal{F}_{i-1}\right] \\
&= A^{d-i}x_i+\sum_{j=1}^{d-i}A^{j-1}Bu_{-j},
\end{aligned}
\end{equation*}
for $i=0,\ldots,d-1$.

The matrices $Z_k$ and $X_k$, $k=N,N-1,\ldots,d$, satisfy the following $\lambda$-dependent Riccati-ZXL equation:
\begin{align}
Z_k
&= A^\top Z_{k+1}A
   + \sigma^2\bar{A}^\top X_{k+1}\bar{A}
   + Q(\lambda) - L_k,
\label{equ:1-6} \\
X_k
&= Z_k+\sum_{i=0}^{d-1}(A^\top)^iL_{k+i}A^i,
\label{equ:1-7}
\end{align}
where
\begin{align}
L_k
&= M_k^\top \Upsilon_k^{-1}M_k,
\label{equ:1-8} \\
\Upsilon_k
&= B^\top Z_{k+1}B
   + \sigma^2\bar{B}^\top X_{k+1}\bar{B}
   + R(\lambda),
\label{equ:1-9} \\
M_k
&= B^\top Z_{k+1}A
   + \sigma^2\bar{B}^\top X_{k+1}\bar{A},
\label{equ:1-10}
\end{align}
with terminal conditions $Z_{N+1}=X_{N+1}=F(\lambda)$. In addition, $L_k$ is set to zero for $k>N$.
\end{theorem}

\subsection{Optimal Parameter Selection Problem}
This subsection addresses the optimal selection of $\lambda$. Based on Lagrange duality \cite{Lu:97}, the constrained stochastic LQR problem is reformulated as a dual maximization problem.

We solve the resulting dual problem via projected gradient ascent. The following theorem characterizes the required dual gradient and the corresponding iteration.
\begin{theorem}\label{thm:1-2}
Consider Problem 1 and suppose that the Slater condition in
Assumption \ref{asp:1-1} holds. For any given $\lambda = (\lambda_1,\ldots,\lambda_r)^\top \in \mathbb{R}_+^r$, let
$\mathbf{u}^*(\lambda)$ be the optimal controller of
\eqref{equ:1-3}. Then
\begin{equation}\label{equ:1-11}
\nabla_\lambda\varphi(\lambda)
=
\left.
\frac{\partial J(\mathbf{u},\lambda)}{\partial\lambda}
\right|_{\mathbf{u}=\mathbf{u}^*(\lambda)} .
\end{equation}
Moreover, the $j$th component of the gradient is given by
\begin{align}
\frac{\partial\varphi(\lambda)}{\partial\lambda_j}
={}&
\mathbb{E}\Bigg[
\sum_{i=0}^{d-1}x_i^\top Q_jx_i
+x_d^\top X_{d,j}x_d
\notag\\
&\quad
-\sum_{i=0}^{d-1}
x_d^\top(A^\top)^iL_{d+i,j}A^i\hat{x}_{d\mid i}
\Bigg]-c_j ,
\label{equ:1-12}
\end{align}
These sensitivity matrices satisfy, for $k=N,N-1,\ldots,d$,
\begin{align}
Z_{k,j}
={}&
A^\top Z_{k+1,j}A
+\sigma^2\bar{A}^\top X_{k+1,j}\bar{A}
+Q_j-L_{k,j},
\label{equ:1-13}\\
X_{k,j}
={}&
Z_{k,j}
+\sum_{i=0}^{d-1}(A^\top)^iL_{k+i,j}A^i,
\label{equ:1-14}\\
L_{k,j}
={}&
M_{k,j}^\top\Upsilon_k^{-1}M_k
+M_k^\top\Upsilon_k^{-1}M_{k,j}
\notag\\
&\quad
-M_k^\top\Upsilon_k^{-1}
\Upsilon_{k,j}\Upsilon_k^{-1}M_k,
\label{equ:1-15}\\
\Upsilon_{k,j}
={}&
B^\top Z_{k+1,j}B
+\sigma^2\bar{B}^\top X_{k+1,j}\bar{B}
+R_j,
\label{equ:1-16}\\
M_{k,j}
={}&
B^\top Z_{k+1,j}A
+\sigma^2\bar{B}^\top X_{k+1,j}\bar{A},
\label{equ:1-17}
\end{align}
with terminal values $Z_{N+1,j}=X_{N+1,j}=F_j$.
Consequently, the projected gradient ascent iteration is
\begin{equation}\label{equ:1-18}
\lambda_{m+1}
=
\bigl[
\lambda_m+\alpha\nabla_\lambda\varphi(\lambda_m)
\bigr]_+,
\end{equation}
where $\alpha>0$ is the step size.
\end{theorem}
\textbf{Proof.} See Appendix A.
\hfill$\Box$

\begin{remark}\label{rem:1-1}
Let $\lambda^*$ be an optimal dual solution, and let $\mathbf{u}^*$ be a corresponding optimal control for Problem~1. Then, by the Karush--Kuhn--Tucker (KKT) conditions, for each $i=1,\ldots,r$, we have
\begin{equation}\label{equ:1-19}
\lambda_i^*\bigl(J_i(\mathbf{u}^*)-c_i\bigr)=0 .
\end{equation}
\end{remark}

The dual function $\varphi(\lambda)$ is concave on $\mathbb R_+^r$. Hence, the dual problem can be solved by the projected gradient ascent algorithm \eqref{equ:1-18}. For the subsequent convergence analysis, the following two lemmas establish the Lipschitz continuity of $\nabla_\lambda\varphi(\lambda)$ on bounded subsets of $\mathbb R_+^r$.

\begin{lemma}\label{lem:1-1}
Let $\Lambda\subset\mathbb R_+^r$ be bounded. Then, for each
$j=1,\ldots,r$, the Riccati-ZXL matrices $Z_k$, $X_k$,
$\Upsilon_k$, $M_k$, $L_k$, and their sensitivity matrices
$Z_{k,j}$, $X_{k,j}$, $\Upsilon_{k,j}$, $M_{k,j}$, $L_{k,j}$,
are uniformly bounded and Lipschitz continuous with respect to
$\lambda$ on $\Lambda$.

In particular, there exist constants $l_j^x>0$ and
$l_{d+i,j}^L>0$ such that, for all $\lambda,\lambda'\in\Lambda$,
\begin{equation*}
\|X_{d,j}(\lambda)-X_{d,j}(\lambda')\|
\le l_j^x\|\lambda-\lambda'\|,
\end{equation*}
and
\begin{equation*}
\|L_{d+i,j}(\lambda)-L_{d+i,j}(\lambda')\|
\le l_{d+i,j}^L\|\lambda-\lambda'\|,
\quad i=0,\ldots,d-1 .
\end{equation*}
\end{lemma}

\textbf{Proof.} See Appendix B.
\hfill$\Box$

\begin{lemma}\label{lem:1-2}
Let $\Lambda\subset\mathbb R_+^r$ be bounded. Then the
gradient $\nabla_\lambda\varphi(\lambda)$ is Lipschitz continuous
on $\Lambda$. Specifically, there exists a constant $l>0$ such that, for any $\lambda,\lambda'\in\Lambda$,
\begin{equation*}
\|\nabla_\lambda\varphi(\lambda)
-\nabla_\lambda\varphi(\lambda')\|
\le
l\|\lambda-\lambda'\|,
\end{equation*}
where
\begin{equation*}
l
=
\sqrt{
\sum_{j=1}^{r}
\left(
l_j^x+
\sum_{i=0}^{d-1}
\|A\|^{2i}l_{d+i,j}^L
\right)^2
}
\,\mathbb E[\|x_d\|^2].
\end{equation*}
\end{lemma}

\textbf{Proof.} See Appendix C.
\hfill$\Box$

Based on \cite{Bo:04}, we obtain the following convergence guarantee.

\begin{corollary}\label{cor:1-1}
Under Assumption \ref{asp:1-1}, suppose that the generated dual sequence remains in a bounded set containing a dual optimal solution $\lambda^*$. Then the projected gradient ascent algorithm \eqref{equ:1-18} with stepsize $0<\alpha\le 1/l$ satisfies $\varphi(\lambda^*)-\varphi(\lambda_n)=O(1/n)$.
\end{corollary}

\section{Infinite-horizon constrained stochastic optimal LQR control}
\subsection{Problem Statement}
In this section, the system we are considering is still described by \eqref{equ:1-1}. In the infinite-horizon case, the cost $(i=0)$ and constraint $(i=1, \cdots, r)$ functions are given by
\begin{equation}\label{equ:1-20}
\tilde{J}_i(\mathbf{u}) = \mathbb{E}\bigl( \sum_{k=0}^{\infty} x_k^\top Q_i x_k + \sum_{k=d}^{\infty} u_{k-d}^\top R_i u_{k-d} \bigr),
\end{equation}
where $\mathbb{E}$ denotes the mathematical expectation with respect to the noise sequence $\{\omega_0, \omega_1, \ldots\}$, and where the weighting matrices $Q_i$ and $R_i$ are as in \eqref{equ:1-2}.

We first introduce the following definitions.
\begin{definition}\label{def:1-1}
The dynamic system \eqref{equ:1-1} is mean-square stabilizable if there exists a feedback controller $u_{k-d} = K \hat{x}_{k|k-d} = K E[x_k |\mathcal{F}_{k-d-1}]$ for $k \geq d$, such that the discrete-time system:
\begin{equation}\label{equ:1-21}
x_{k+1} = \big(A + \omega_k \bar{A}\big) x_k + \big(B + \omega_k \bar{B}\big) K \hat{x}_{k \mid k-d}
\end{equation}
is asymptotically mean-square stable. In other words, for any initial values $x_0$ and $u_{k-d}$, $k = 0, \ldots, d-1$, the state satisfies $\lim\limits_{k \rightarrow \infty} \mathbb{E}(x_k^\top x_k) = 0$.
\end{definition}

\begin{definition}\label{def:1-2}
The stochastic system
\begin{equation}\label{equ:1-22}
x_{k+1} = (A + \omega_k \bar{A}) x_k, \quad y_k = C x_k
\end{equation}
is said to be exactly observable (or, for short, $(A, \bar{A}, C)$ is said to be exactly observable)
if for any $N \ge n$,
\begin{equation*}
y_k \equiv 0,\ \text{a.s.}\ \forall\, 0 \le k \le N \ \Rightarrow\ x_0 = 0.
\end{equation*}
\end{definition}

Given $c_1, \dots, c_r \in \mathbb{R}$, the constrained stochastic optimal LQR control problem for the infinite-horizon case as follows:

\noindent\textbf{Problem 2}
Find an $\mathcal{F}_{k-1}$-adapted $u_k$ such that $\tilde{J}_0(\mathbf{u})$ is minimized under the quadratic constraints $\tilde{J}_i(\mathbf{u}) \leq c_i(i=1, \ldots, r)$, subject to the dynamical system \eqref{equ:1-1} being mean-square stabilizable.

\subsection{Infinite-Horizon Stochastic Optimal Control}
\textbf{Problem 2} is an infinite-dimensional optimization problem. Nevertheless, under the Slater condition and our assumptions, strong duality still holds, thereby enabling the determination of optimal controls. The Slater condition is defined as follows:
\begin{assumption}[Slater Condition]\label{asp:1-2}
There exists an $\mathcal{F}_{k-1}$-adapted $\bar{u}_k$ such that for all $i = 1, \ldots, r$,
\begin{equation*}
\tilde{J}_i(\bar{\mathbf{u}}) < c_i.
\end{equation*}
\end{assumption}

For $\lambda = (\lambda_1,\ldots,\lambda_r)^\top \in \mathbb{R}_+^r$, define the Lagrangian:
\begin{equation*}
\tilde{J}(\mathbf{u},\lambda)
=
\tilde{J}_0(\mathbf{u})
+
\sum_{i=1}^{r}\lambda_i\bigl(\tilde{J}_i(\mathbf{u})-c_i\bigr),
\end{equation*}
where $\lambda$ is referred to as the Lagrange multiplier. It follows that
\begin{equation*}
\begin{aligned}
\tilde{J}(\mathbf{u},\lambda)
=
\mathbb{E}\Bigg[
&\sum_{k=0}^{\infty} x_k^\top Q(\lambda)x_k
+
\sum_{k=d}^{\infty} u_{k-d}^\top R(\lambda)u_{k-d}
\Bigg]
-\lambda^\top c,
\end{aligned}
\end{equation*}
where $c=\left(c_1,\ldots,c_r\right)^\top$, and $Q(\lambda) = Q_0 + \sum\limits_{i=1}^{r} \lambda_i Q_i$, with $R(\lambda)$ defined analogously.

To guarantee the solvability of \textbf{Problem 2}, we impose the following condition. For every $\lambda\in\mathbb{R}_+^r$, it follows from
$R_0\succ0$, $R_i\succeq0$ and $Q_i\succeq0$ that $R(\lambda)\succ0$ and $Q(\lambda)\succeq0$. Hence, there exists a matrix $C_\lambda$ such that $Q(\lambda)=C_\lambda^\top C_\lambda$, which ensures the uniqueness of the optimal controller. The following assumption is standard for mean-square stabilization.

\begin{assumption}\label{asp:1-3}
For every $\lambda\in\mathbb{R}_+^r$, $(A,\bar A,C_\lambda)$ is exactly observable.
\end{assumption}

In order to reflect the finite horizon $N$ in the finite-horizon
stochastic LQR problem, we denote the matrices $Z_k$ and $X_k$ in the
Riccati--ZXL equation \eqref{equ:1-6}--\eqref{equ:1-10} by
$Z_k(N)$ and $X_k(N)$, respectively. For each fixed
$\lambda\in\mathbb{R}_+^r$, these matrices depend on $\lambda$
through $Q(\lambda)$ and $R(\lambda)$. According to \cite{Zh:15}, the following lemma holds.

\begin{lemma}\label{lem:1-3}
For any fixed $\lambda\in\mathbb{R}_+^r$ and any $k \geq 0$,
$Z_k(N)$ and $X_k(N)$ are convergent as $N \rightarrow \infty$, i.e.,
\begin{equation*}
Z \doteq \lim_{N \to \infty} Z_k(N),
\quad
X \doteq \lim_{N \to \infty} X_k(N)
\end{equation*}
exist and are independent of $k$. Moreover, $Z$ and $X$ satisfy the
following Riccati--ZXL equation:
\begin{align}
Z
&=
A^\top Z A
+
\sigma^2 \bar{A}^\top X \bar{A}
+
Q(\lambda)
-
L,
\label{equ:1-23}
\\
X
&=
Z
+
\sum_{i=0}^{d-1}
\bigl(A^\top\bigr)^i
L
A^i,
\label{equ:1-24}
\end{align}
with
\begin{align}
L
&=
M^\top
\Upsilon^{-1}
M,
\label{equ:1-25}
\\
\Upsilon
&=
B^\top Z B
+
\sigma^2 \bar{B}^\top X\bar{B}
+
R(\lambda),
\label{equ:1-26}
\\
M
&=
B^\top Z A
+
\sigma^2 \bar{B}^\top X\bar{A}.
\label{equ:1-27}
\end{align}
\end{lemma}

We also consider the Lagrange dual problem. Given any $\lambda\in\mathbb{R}_+^r$, we have the associated dual function
\begin{equation}\label{equ:1-28}
\psi(\lambda)=\min _\mathbf{u} \tilde{J}(\mathbf{u}, \lambda)
\end{equation}

This is a parameterized, unconstrained infinite-horizon stochastic optimal LQR control problem with delay and multiplicative noise. The following theorem is a direct result from \cite{Zh:15}.
\begin{theorem}\label{thm:1-3}
Under Assumption \ref{asp:1-3}, for every
$\lambda\in\mathbb{R}_+^r$, the system \eqref{equ:1-1} is stabilizable
in the mean-square sense if and only if there exists a unique solution
to \eqref{equ:1-23}--\eqref{equ:1-27} satisfying $Z>0$. In this case,
the optimal $\mathcal{F}_{k-1}$-adapted controller for $\psi(\lambda)$ is
\begin{equation}\label{equ:1-29}
u_k=-\Upsilon^{-1}M\hat{x}_{k+d\mid k},\quad k\geq 0 .
\end{equation}
The associated optimal cost can be expressed as
\begin{equation}\label{equ:1-30}
\begin{aligned}
\psi(\lambda)
={}& x_0^\top Zx_0
-\sum_{k=0}^{d-1}u_{k-d}^\top R(\lambda)u_{k-d}
-\lambda^\top c                                      \\
&+\sum_{k=0}^{d-1}\mathbb{E}\Big[
\big(u_{k-d}+\Upsilon^{-1}M\hat{x}_{k\mid k-d}\big)^\top \\
&\hspace{3.6em}\times
\Upsilon
\big(u_{k-d}+\Upsilon^{-1}M\hat{x}_{k\mid k-d}\big)
\Big],
\end{aligned}
\end{equation}
where
\begin{equation*}
\hat{x}_{k\mid k-d}
=
A^k x_0+\sum_{j=0}^{k-1}A^{k-1-j}B u_{j-d},
\end{equation*}
for $k=0,\ldots,d-1$.
\end{theorem}

\subsection{Optimal Parameter Selection Problem}
In this subsection, we formulate the optimal parameter selection problem via duality theory and, similarly, develop a projected gradient ascent method for solving the resulting dual problem.

For the infinite-horizon problem, the minimization with respect to $\mathbf u$ is understood to be taken over the class of admissible mean-square stabilizing adapted controls. More precisely, an admissible control is an adapted control sequence for which the resulting closed-loop state process is mean-square stable and the performance functionals in \eqref{equ:1-20} are well defined and finite. Thus, the mean-square stabilizing requirement is incorporated into the admissible control class, rather than being treated as an inequality constraint to be dualized. Accordingly, Lagrange multipliers are introduced only for the constraints $\tilde J_i(\mathbf u)\le c_i$, $i=1,\ldots,r$.

Although the infinite-horizon problem is infinite-dimensional, it falls
within the standard convex-analytic framework for optimization over a
function space. Indeed, since the state depends linearly on the control
sequence and $Q_i\succeq0$, $R_i\succeq0$, the functionals $\tilde J_i(\mathbf u)$, $i=0,1,\ldots,r$, are convex quadratic functionals of $\mathbf u$. Moreover, the constraints $\tilde J_i(\mathbf u)\le c_i$, $i=1,\ldots,r$, are finitely many scalar convex inequality constraints. Hence, by the standard strong duality result for infinite-dimensional convex programs, the Slater condition in Assumption \ref{asp:1-2} implies zero duality gap.

\begin{proposition}\label{prop:1-1}
Under Assumption \ref{asp:1-2}, the following strong duality relation holds:
\begin{equation}\label{equ:1-31}
\begin{aligned}
\min_{\mathbf u}
\max_{\lambda\ge0}
\tilde J(\mathbf u,\lambda)
&=
\max_{\lambda\ge0}
\min_{\mathbf u}
\tilde J(\mathbf u,\lambda)  \\
&=
\max_{\lambda\ge0}\psi(\lambda).
\end{aligned}
\end{equation}
\end{proposition}

Based on Proposition \ref{prop:1-1}, we obtain the following result.
\begin{theorem}\label{thm:1-4}
Consider Problem 2 and suppose that Assumptions
\ref{asp:1-2} and \ref{asp:1-3} hold. For any given
$\lambda=(\lambda_1,\ldots,\lambda_r)^\top\in\mathbb{R}_+^r$, let
$\mathbf{u}^*(\lambda)$ be the optimal controller of
\eqref{equ:1-28}. Then
\begin{equation}\label{equ:1-32}
\nabla_\lambda\psi(\lambda)
=
\left.
\frac{\partial \tilde{J}(\mathbf{u},\lambda)}
{\partial\lambda}
\right|_{\mathbf{u}=\mathbf{u}^*(\lambda)} .
\end{equation}
Moreover, the $j$th component of the gradient is given by
\begin{align}
\frac{\partial\psi(\lambda)}{\partial\lambda_j}
={}&
x_0^\top Z_jx_0-c_j
\notag\\
&+
\sum_{k=0}^{d-1}\mathbb{E}\Big[
2u_{k-d}^\top M_j\hat{x}_{k\mid k-d}
\notag\\
&\quad
+\hat{x}_{k\mid k-d}^\top
L_j\hat{x}_{k\mid k-d}
\notag\\
&\quad
+u_{k-d}^\top
(\Upsilon_j-R_j)u_{k-d}
\Big].
\label{equ:1-33}
\end{align}
These sensitivity matrices satisfy
\begin{align}
Z_j
={}&
A^\top Z_jA
+\sigma^2\bar{A}^\top X_j\bar{A}
+Q_j-L_j,
\label{equ:1-34}
\\
X_j
={}&
Z_j
+\sum_{i=0}^{d-1}(A^\top)^iL_jA^i,
\label{equ:1-35}
\\
L_j
={}&
M_j^\top\Upsilon^{-1}M
+M^\top\Upsilon^{-1}M_j
\notag\\
&-
M^\top\Upsilon^{-1}
\Upsilon_j\Upsilon^{-1}M,
\label{equ:1-36}
\\
\Upsilon_j
={}&
B^\top Z_jB
+\sigma^2\bar{B}^\top X_j\bar{B}
+R_j,
\label{equ:1-37}
\\
M_j
={}&
B^\top Z_jA
+\sigma^2\bar{B}^\top X_j\bar{A}.
\label{equ:1-38}
\end{align}
Consequently, the projected gradient ascent iteration for the dual
variable is
\begin{equation}\label{equ:1-39}
\lambda_{m+1}
=
\bigl[
\lambda_m+\alpha\nabla_\lambda\psi(\lambda_m)
\bigr]_+,
\end{equation}
where $\alpha>0$ is the step size.
\end{theorem}

\textbf{Proof.} See Appendix D.
\hfill$\Box$

\begin{remark}\label{rem:1-2}
Let $\lambda^*$ be an optimal dual solution, and let
$\mathbf{u}^*$ be the corresponding optimal controller for
Problem~2. Then, by the Karush--Kuhn--Tucker (KKT)
conditions, for each $i=1,\ldots,r$, we have
\begin{equation}\label{equ:1-40}
\lambda_i^*
\bigl(\tilde{J}_i(\mathbf{u}^*)-c_i\bigr)
=0 .
\end{equation}
\end{remark}

For the subsequent convergence analysis, the following two lemmas establish the Lipschitz continuity of $\nabla_\lambda\psi(\lambda)$ on bounded subsets of $\mathbb R_+^r$.

\begin{lemma}\label{lem:1-4}
Let $\Lambda\subset\mathbb R_+^r$ be bounded. Then, for each
$j=1,\ldots,r$, the matrix functions $Z_j$, $M_j$, $L_j$, and
$\Upsilon_j-R_j$ are uniformly bounded and Lipschitz continuous with
respect to $\lambda$ on $\Lambda$.

In particular, there exist constants $l_j^Z>0$, $l_j^M>0$,
$l_j^L>0$, and $l_j^\Upsilon>0$ such that, for all
$\lambda,\lambda'\in\Lambda$,
\begin{align*}
\|Z_j(\lambda)-Z_j(\lambda')\|
&\le l_j^Z\|\lambda-\lambda'\|,\\
\|M_j(\lambda)-M_j(\lambda')\|
&\le l_j^M\|\lambda-\lambda'\|,\\
\|L_j(\lambda)-L_j(\lambda')\|
&\le l_j^L\|\lambda-\lambda'\|,\\
\|[\Upsilon_j(\lambda)-R_j]
-[\Upsilon_j(\lambda')-R_j]\|
&\le l_j^\Upsilon\|\lambda-\lambda'\|.
\end{align*}
\end{lemma}

\textbf{Proof.} See Appendix E.
\hfill$\Box$

\begin{lemma}\label{lem:1-5}
Let $\Lambda\subset\mathbb R_+^r$ be bounded. Then
$\nabla_\lambda\psi(\lambda)$ is Lipschitz continuous on
$\Lambda$. Specifically, there exists a constant $l'>0$,
depending possibly on $\Lambda$, such that, for any
$\lambda,\lambda'\in\Lambda$,
\begin{equation*}
\|
\nabla_\lambda\psi(\lambda)
-
\nabla_\lambda\psi(\lambda')
\|
\le
l'\|\lambda-\lambda'\|.
\end{equation*}
Moreover, one may take
\begin{equation*}
l'
=
\sqrt{
\sum_{j=1}^{r}(l_j')^2
},
\end{equation*}
where
\begin{equation*}
l_j'
=
l_j^Z\|x_0\|^2
+
\sum_{k=0}^{d-1}
\left(
2l_j^M\mathbb E_k^{(1)}
+
l_j^L\mathbb E_k^{(2)}
+
l_j^\Upsilon\mathbb E_k^{(3)}
\right),
\end{equation*}
and the following finite constants are determined by the
initial data:
\begin{equation*}
\begin{cases}
\mathbb E_k^{(1)}
=
\mathbb E\bigl[
\|u_{k-d}\|\,
\|\hat x_{k\mid k-d}\|
\bigr],\\[2mm]
\mathbb E_k^{(2)}
=
\mathbb E\bigl[
\|\hat x_{k\mid k-d}\|^2
\bigr],\\[2mm]
\mathbb E_k^{(3)}
=
\mathbb E\bigl[
\|u_{k-d}\|^2
\bigr].
\end{cases}
\end{equation*}
\end{lemma}

\textbf{Proof.} See Appendix F.
\hfill$\Box$

Based on \cite{Bo:04}, we obtain the following convergence guarantee.

\begin{corollary}\label{cor:1-2}
Under Assumptions \ref{asp:1-2}--\ref{asp:1-3}, suppose that the generated dual sequence remains in a bounded set containing a dual optimal solution $\lambda^*$. Then the projected gradient ascent algorithm \eqref{equ:1-33} with stepsize $0<\alpha\le 1/l'$ satisfies $\psi(\lambda^*)-\psi(\lambda_n)=O(1/n)$.
\end{corollary}

\section{Numerical Examples}
\subsection{The Finite-Horizon Case}

Consider the following discrete-time system:
\begin{equation}\label{equ:1-41}
x_{k+1}=\left(1+\omega_k\right)x_k+\left(2+2\omega_k\right)u_{k-1},
\end{equation}
where $\sigma^2=1$, with the initial values $x_0=1$ and $u_{-1}=-1$. Problem~1 is formulated as follows:
\begin{equation*}
\left\{
\begin{aligned}
\min_{\mathbf{u}}\quad
& J_0(\mathbf{u}) =
\mathbb{E}\left[
\sum_{k=0}^{2}2\|x_k\|^2
+\sum_{k=1}^{2}5\|u_{k-1}\|^2
+5\|x_3\|^2
\right] \\
\text{s.t.}\quad
& J_1(\mathbf{u}) =
\mathbb{E}\left[
\sum_{k=0}^{2}2\|x_k\|^2
+\sum_{k=1}^{2}3\|u_{k-1}\|^2
+\|x_3\|^2
\right] \\
& \hspace{18mm} \leq 13.25,\\
& x_k,\ \mathcal{F}_{k-1}\text{-adapted }u_k
\text{ satisfy } \eqref{equ:1-41}.
\end{aligned}
\right.
\end{equation*}

We choose the learning rate $\alpha=0.01$, the initial value $\lambda_0=0$, and the error tolerance $e=10^{-9}$.

By the projected gradient ascent algorithm \eqref{equ:1-18}, we obtain $\lambda^*=2.2313$. With this multiplier, the Riccati-ZXL equation is solved. Then, by \eqref{equ:1-4}, the optimal controller for Problem 1 is obtained as
\begin{equation*}
u_0^*(\lambda^*)=-0.4554\hat{x}_{1\mid 0},\quad
u_1^*(\lambda^*)=-0.4159\hat{x}_{2\mid 1}.
\end{equation*}

The corresponding optimal value of the cost \eqref{equ:1-5} is
$J_0(\mathbf{u}^*(\lambda^*))=22.30$, and the constraint cost is
$J_1(\mathbf{u}^*(\lambda^*))=13.25$. Thus, the KKT condition
\eqref{equ:1-19} is satisfied.

\subsection{The Infinite-Horizon Case}
Consider the following discrete-time system:
\begin{equation}\label{equ:1-42}
x_{k+1}
=\left(1.3+0.1\omega_k\right)x_k
+\left(0.2+0.1\omega_k\right)u_{k-1},
\end{equation}
where $\sigma^2=1$, with the initial values $x_0=1$ and
$u_{-1}=-1$. Problem~2 is formulated as follows:
\begin{equation*}
\left\{
\begin{aligned}
\min_{\mathbf{u}}\quad
& \tilde{J}_0(\mathbf{u})
=
\mathbb{E}\left[
\sum_{k=0}^{\infty}\|x_k\|^2
+\sum_{k=1}^{\infty}\|u_{k-1}\|^2
\right] \\
\text{s.t.}\quad
& \tilde{J}_1(\mathbf{u})
=
\mathbb{E}\left[
\sum_{k=0}^{\infty}0.5\|x_k\|^2
+\sum_{k=1}^{\infty}2\|u_{k-1}\|^2
\right] \\
& \hspace{18mm} \leq 49.35,\\
& \tilde{J}_2(\mathbf{u})
=
\mathbb{E}\left[
\sum_{k=0}^{\infty}0.1\|x_k\|^2
+\sum_{k=1}^{\infty}1.9\|u_{k-1}\|^2
\right] \\
& \hspace{18mm} \leq 45.21,\\
& x_k,\ \mathcal{F}_{k-1}\text{-adapted }u_k
\text{ satisfy } \eqref{equ:1-42},\\
& \text{and the closed-loop system is mean-square stable.}
\end{aligned}
\right.
\end{equation*}

We choose the learning rate $\alpha=0.001$, the initial value
$\lambda_0=(0,0)^\top$, and the error tolerance $e=10^{-9}$.

By the projected gradient ascent algorithm \eqref{equ:1-39}, we obtain
$\lambda^*=(0.1712,0.3178)^\top$. By plugging it into the Riccati-ZXL equation \eqref{equ:1-23}--\eqref{equ:1-27}, we obtain
$$Z=41.0826,\quad X=71.2596.$$

It is assured that a unique optimal controller exists, according to Theorem \ref{thm:1-3}, which stabilizes system \eqref{equ:1-42} in the mean-square sense. Then, by \eqref{equ:1-29}, the optimal controller for Problem 2 is obtained as
\begin{equation*}
u_k^*(\lambda^*)=-2.6485\hat{x}_{k+1\mid k},\quad k\geq0 .
\end{equation*}

\begin{figure}[t]
\centering
\includegraphics[width=3.7in]{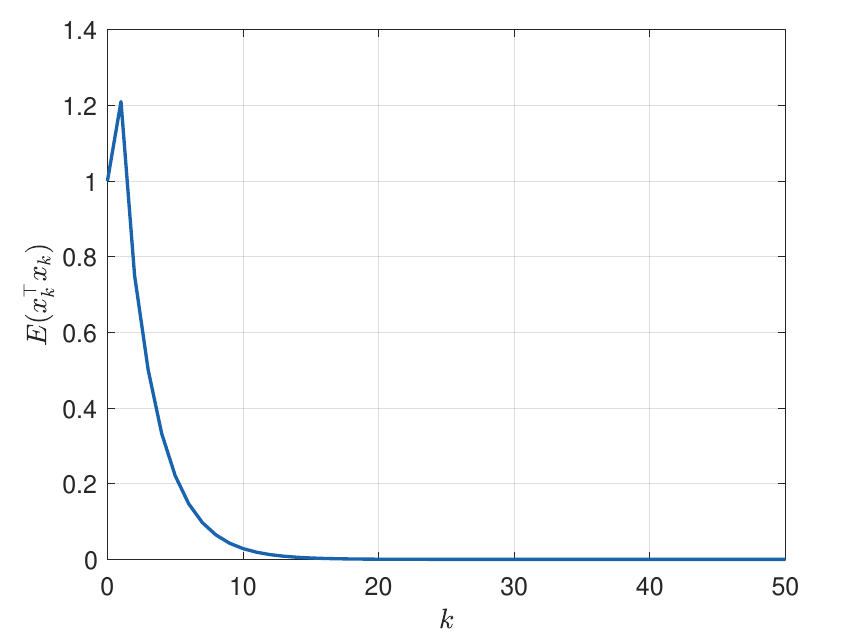}
\caption{Dynamic Behavior of $E\left(x_k^T x_k\right)$.}
\label{untitled1}
\end{figure}

The corresponding optimal value of the cost \eqref{equ:1-30} is
$\tilde{J}_0(\mathbf{u}^*(\lambda^*)) = 28.01$, with constraint costs
$\tilde{J}_1(\mathbf{u}^*(\lambda^*)) = 49.35$ and
$\tilde{J}_2(\mathbf{u}^*(\lambda^*)) = 45.21$. Thus, the KKT condition \eqref{equ:1-40} is satisfied. Fig. \ref{untitled1} presents a simulation outcome for the designed controller, demonstrating that the regulated state attains asymptotic mean-square stability.

\section{Conclusion}
In this article, we study stochastic optimal LQR control for discrete-time systems with delays and quadratic constraints. Using duality theory, we leverage strong duality to transform the primal problem into its dual problem, which consists of a parameterized unconstrained stochastic LQR problem together with an optimal parameter selection problem. By solving the Riccati-ZXL equation, we obtain the parameterized optimal controller and its corresponding cost. The optimal parameters are determined via a projected gradient ascent algorithm, which achieves an $O(1/n)$ dual objective gap on any bounded set containing a dual optimal solution. Future work will address stochastic optimal control for partially observed systems with delay and quadratic constraints.

\section*{Appendix A}
\section*{Proof of theorem 2}
Under the conditions stated in the theorem, Problem 1 is convex in $\mathbf{u}$. By the Slater condition in
Assumption \ref{asp:1-1}, strong duality holds. For any fixed
$\lambda\geq0$, Theorem \ref{thm:1-1} gives the unique optimal
controller $\mathbf{u}^*(\lambda)$ and the corresponding value
\begin{equation*}
\varphi(\lambda)=\min_{\mathbf{u}}J(\mathbf{u},\lambda).
\end{equation*}

Since $J(\mathbf{u},\lambda)$ is differentiable with respect to
$\lambda$, Danskin's theorem implies
\begin{equation*}
\nabla_\lambda\varphi(\lambda)
=
\left.
\frac{\partial J(\mathbf{u},\lambda)}{\partial\lambda}
\right|_{\mathbf{u}=\mathbf{u}^*(\lambda)} ,
\end{equation*}
which proves \eqref{equ:1-11}. Differentiating the optimal
cost in Theorem \ref{thm:1-1} with respect to $\lambda_j$
yields \eqref{equ:1-12}. Similarly, differentiating the
backward recursions for $Z_k$, $X_k$, $L_k$, $\Upsilon_k$,
and $M_k$ gives \eqref{equ:1-13}--\eqref{equ:1-17}, with
terminal values $Z_{N+1,j}=X_{N+1,j}=F_j$. The projected
gradient ascent update \eqref{equ:1-18} then follows directly.
This completes the proof.

\section*{Appendix B}
\section*{Proof of Lemma~\ref{lem:1-1}}

Let $\Lambda\subset\mathbb R_+^r$ be bounded. Since
$Q(\lambda)$, $R(\lambda)$, and $F(\lambda)$ are affine in
$\lambda$, they are uniformly bounded and Lipschitz continuous on
$\Lambda$. We prove the assertion by backward induction.

At $k=N+1$, the terminal conditions
$Z_{N+1}=X_{N+1}=F(\lambda)$ and
$Z_{N+1,j}=X_{N+1,j}=F_j$ imply the desired properties.

Assume that, for some $t\in\{d,\ldots,N\}$, the matrices
$Z_k$, $X_k$, $\Upsilon_k$, $M_k$, $L_k$, and their sensitivity
matrices $Z_{k,j}$, $X_{k,j}$, $\Upsilon_{k,j}$, $M_{k,j}$,
$L_{k,j}$ are uniformly bounded and Lipschitz continuous on
$\Lambda$ for all admissible $k\ge t+1$. Then \eqref{equ:1-9}
implies that $\Upsilon_t$ is uniformly bounded and Lipschitz
continuous on $\Lambda$.

Moreover, since $Q(\lambda)\succeq0$, $F(\lambda)\succeq0$, and
$R(\lambda)\succeq R_0\succ0$ for all $\lambda\in\mathbb R_+^r$,
it follows from \cite{Zh:15} that
$Z_{t+1}(\lambda)\succeq0$ and $X_{t+1}(\lambda)\succeq0$.
Hence $\Upsilon_t(\lambda)\succeq R_0\succ0$ on $\Lambda$, and
therefore $\Upsilon_t^{-1}$ is uniformly bounded.

For any $\lambda,\lambda'\in\Lambda$,
\begin{equation*}
\begin{aligned}
&\Upsilon_t^{-1}(\lambda)-\Upsilon_t^{-1}(\lambda')  \\
&\quad =
\Upsilon_t^{-1}(\lambda)
\bigl(\Upsilon_t(\lambda')-\Upsilon_t(\lambda)\bigr)
\Upsilon_t^{-1}(\lambda') .
\end{aligned}
\end{equation*}
Thus, using the uniform boundedness of $\Upsilon_t^{-1}$ and the
Lipschitz continuity of $\Upsilon_t$, we obtain
\begin{equation*}
\|\Upsilon_t^{-1}(\lambda)-\Upsilon_t^{-1}(\lambda')\|
\le C\|\lambda-\lambda'\|
\end{equation*}
for some constant $C>0$. Hence $\Upsilon_t^{-1}$ is Lipschitz
continuous on $\Lambda$.

By \eqref{equ:1-10}, $M_t$ is uniformly bounded and Lipschitz
continuous on $\Lambda$. Since
$L_t=M_t^\top\Upsilon_t^{-1}M_t$, the same property holds for
$L_t$. Then \eqref{equ:1-6} and \eqref{equ:1-7} imply that
$Z_t$ and $X_t$ are uniformly bounded and Lipschitz continuous on
$\Lambda$.

It remains to consider the sensitivity matrices. From
\eqref{equ:1-16} and \eqref{equ:1-17}, together with the induction
hypothesis and the fact that $R_j$ is constant, $\Upsilon_{t,j}$
and $M_{t,j}$ are uniformly bounded and Lipschitz continuous on
$\Lambda$. Furthermore, by \eqref{equ:1-15},
\begin{equation*}
\begin{aligned}
L_{t,j}
={}&
M_{t,j}^\top\Upsilon_t^{-1}M_t
+M_t^\top\Upsilon_t^{-1}M_{t,j}  \\
&-
M_t^\top\Upsilon_t^{-1}
\Upsilon_{t,j}\Upsilon_t^{-1}M_t .
\end{aligned}
\end{equation*}
All factors on the right-hand side are uniformly bounded and
Lipschitz continuous on $\Lambda$; hence so is $L_{t,j}$. Finally,
\eqref{equ:1-13} and \eqref{equ:1-14} yield the same properties
for $Z_{t,j}$ and $X_{t,j}$.

The induction is complete. Therefore, the Riccati-ZXL matrices
$Z_k$, $X_k$, $\Upsilon_k$, $M_k$, $L_k$, and their sensitivity
matrices $Z_{k,j}$, $X_{k,j}$, $\Upsilon_{k,j}$, $M_{k,j}$,
$L_{k,j}$ are uniformly bounded and Lipschitz continuous with
respect to $\lambda$ on $\Lambda$.

In particular, there exist constants $l_j^x>0$ and
$l_{d+i,j}^L>0$, $i=0,\ldots,d-1$, such that, for all
$\lambda,\lambda'\in\Lambda$,
\begin{equation*}
\|X_{d,j}(\lambda)-X_{d,j}(\lambda')\|
\le l_j^x\|\lambda-\lambda'\|,
\end{equation*}
and
\begin{equation*}
\|L_{d+i,j}(\lambda)-L_{d+i,j}(\lambda')\|
\le l_{d+i,j}^L\|\lambda-\lambda'\|,
\quad i=0,\ldots,d-1 .
\end{equation*}
This completes the proof.

\section*{Appendix C}
\section*{Proof of Lemma \ref{lem:1-2}}

Let $\lambda,\lambda'\in\Lambda$. Substituting $\lambda$ and
$\lambda'$ into \eqref{equ:1-12}, and using the triangle
inequality, the Cauchy--Schwarz inequality, and standard norm
bounds, we obtain, for each $j=1,\ldots,r$,
\begin{equation*}
\begin{aligned}
&\left|
\frac{\partial\varphi}{\partial\lambda_j}(\lambda)
-
\frac{\partial\varphi}{\partial\lambda_j}(\lambda')
\right|  \\
&\le
\|X_{d,j}(\lambda)-X_{d,j}(\lambda')\|
\,\mathbb E[\|x_d\|^2]  \\
&\quad+
\sum_{i=0}^{d-1}
\left\|
(A^\top)^i
\bigl(L_{d+i,j}(\lambda)-L_{d+i,j}(\lambda')\bigr)
A^i
\right\|  \\
&\qquad\qquad\times
\mathbb E[\|x_d\|\,\|\hat x_{d|i}\|].
\end{aligned}
\end{equation*}

Since
\begin{equation*}
\hat x_{d|i}=\mathbb E[x_d\mid\mathcal F_{i-1}],
\end{equation*}
Jensen's inequality and the Cauchy--Schwarz inequality imply
\begin{equation*}
\mathbb E[\|x_d\|\,\|\hat x_{d|i}\|]
\le
\mathbb E[\|x_d\|^2].
\end{equation*}

By Lemma~\ref{lem:1-1}, there exist constants $l_j^x>0$ and
$l_{d+i,j}^L>0$, $i=0,\ldots,d-1$, such that
\begin{equation*}
\begin{aligned}
\|X_{d,j}(\lambda)-X_{d,j}(\lambda')\|
&\le l_j^x\|\lambda-\lambda'\|,\\
\|L_{d+i,j}(\lambda)-L_{d+i,j}(\lambda')\|
&\le l_{d+i,j}^L\|\lambda-\lambda'\|.
\end{aligned}
\end{equation*}
Moreover,
\begin{equation*}
\begin{aligned}
&\left\|
(A^\top)^i
\bigl(L_{d+i,j}(\lambda)-L_{d+i,j}(\lambda')\bigr)
A^i
\right\|  \\
&\quad\le
\|A\|^{2i}l_{d+i,j}^L
\|\lambda-\lambda'\|.
\end{aligned}
\end{equation*}

Combining the above estimates gives
\begin{equation*}
\begin{aligned}
&\left|
\frac{\partial\varphi}{\partial\lambda_j}(\lambda)
-
\frac{\partial\varphi}{\partial\lambda_j}(\lambda')
\right|  \\
&\le
\left(
l_j^x+
\sum_{i=0}^{d-1}
\|A\|^{2i}l_{d+i,j}^L
\right)
\mathbb E[\|x_d\|^2]\,
\|\lambda-\lambda'\|.
\end{aligned}
\end{equation*}
Therefore,
\begin{equation*}
\begin{aligned}
&\|\nabla_\lambda\varphi(\lambda)
-\nabla_\lambda\varphi(\lambda')\|  \\
&\le
\sqrt{
\sum_{j=1}^{r}
\left(
l_j^x+
\sum_{i=0}^{d-1}
\|A\|^{2i}l_{d+i,j}^L
\right)^2
}
\,\mathbb E[\|x_d\|^2]\,
\|\lambda-\lambda'\|.
\end{aligned}
\end{equation*}
Hence, $\nabla_\lambda\varphi(\lambda)$ is Lipschitz continuous
on $\Lambda$ with Lipschitz constant
\begin{equation*}
\begin{aligned}
l
={}&
\sqrt{
\sum_{j=1}^{r}
\left(
l_j^x+
\sum_{i=0}^{d-1}
\|A\|^{2i}l_{d+i,j}^L
\right)^2
}  \\
&\times \mathbb E[\|x_d\|^2].
\end{aligned}
\end{equation*}
This completes the proof.

\section*{Appendix D}
\section*{Proof of Theorem \ref{thm:1-4}}

The proof follows from the same argument as that of Theorem
\ref{thm:1-2}, and is therefore omitted.

\section*{Appendix E}
\section*{Proof of Lemma~\ref{lem:1-4}}

By Lemmas \ref{lem:1-1} and \ref{lem:1-3}, the limits of the
finite-horizon Riccati-ZXL sensitivity matrices are uniformly
bounded and Lipschitz continuous on any bounded set
$\Lambda\subset\mathbb R_+^r$. Hence the algebraic sensitivity
matrices appearing in \eqref{equ:1-33},
\[
Z_j,\quad M_j,\quad L_j,\quad \Upsilon_j-R_j,
\]
are uniformly bounded and Lipschitz continuous on $\Lambda$.
This completes the proof.

\section*{Appendix F}
\section*{Proof of Lemma \ref{lem:1-5}}

The proof is similar to that of Lemma \ref{lem:1-2}. Applying
the gradient formula \eqref{equ:1-33} and the Lipschitz
estimates in Lemma \ref{lem:1-4}, one obtains, for each
$j=1,\ldots,r$,
\[
\left|
\frac{\partial\psi(\lambda)}{\partial\lambda_j}
-
\frac{\partial\psi(\lambda')}{\partial\lambda_j}
\right|
\le
l_j'\|\lambda-\lambda'\|.
\]
Taking the Euclidean norm over all components gives
\[
\|
\nabla_\lambda\psi(\lambda)
-
\nabla_\lambda\psi(\lambda')
\|
\le
\left(
\sum_{j=1}^r(l_j')^2
\right)^{1/2}
\|\lambda-\lambda'\|,
\]
which proves the desired result.

\end{document}